\documentclass[11pt]{article}

\usepackage[T1]{fontenc}
\usepackage[utf8]{inputenc}
\usepackage{amsmath,amssymb,amsthm}
\usepackage{enumitem}
\usepackage{microtype}
\usepackage[a4paper,margin=2.7cm]{geometry}
\usepackage{hyperref}
\usepackage{authblk}

\newtheorem{theorem}{Theorem}[section]
\newtheorem{proposition}[theorem]{Proposition}
\newtheorem{lemma}[theorem]{Lemma}
\newtheorem{assumption}[theorem]{Assumption}
\theoremstyle{definition}
\newtheorem{remark}[theorem]{Remark}
\newtheorem{definition}[theorem]{Definition}

\numberwithin{equation}{section}

\newcommand{\R}{\mathbb{R}}
\newcommand{\rank}{\operatorname{rank}}
\newcommand{\im}{\operatorname{im}}

\title{Low-rank iterative computation of available storage for stochastic passive systems}
\author[1]{Julia Ackermann}
\author[2]{Timo Reis}

\affil[1]{Department of Mathematics and Natural Sciences,
University of Wuppertal,
Gaußstraße 20, 42119 Wuppertal, Germany\\
\texttt{jackermann@uni-wuppertal.de}
}
\affil[2]{Systems Theory and Partial Differential Equations Group,
Institute of Mathematics, Technische Universität Ilmenau,
98693 Ilmenau, Germany\\
\texttt{timo.reis@tu-ilmenau.de}}
\date{}
\begin{document}
\maketitle

\begin{abstract}
We develop a sparse low-rank method for computing available storage of
finite-dimensional stochastic passive systems with multiplicative noise.
Available storage is the Loewner-minimal positive semidefinite solution of a
stochastic KYP inequality; singular feedthrough prevents a direct Riccati
reduction, while control-dependent diffusion leads to a rational Riccati
equation. A passivity-based loop transformation preserves the storage solutions
and provides stable dynamics. Structural deflation then removes singular input
directions, preserves the minimal solution, and transports a stabilizing
feedback. Under terminal regular feasibility, Newton's method is initialized by
one stochastic Lyapunov solve. For large sparse problems, the extracted Riccati
equation is realized with projected operators, signed low-rank factors, and
low-rank generalized-Lyapunov solves. The method is illustrated numerically on a sparse stochastic benchmark problem of large dimension.
\end{abstract}


\section{Introduction}

Available storage is the maximal amount of energy that can be extracted from a
system starting from a prescribed initial state. We consider the stochastic
linear input-state-output system
\begin{equation}
\label{eq:system}
\begin{aligned}
dx(t)&=(Ax(t)+Bu(t))\,dt
 +\sum_{j=1}^q(A_jx(t)+B_ju(t))\,dW_j(t),\\
y(t)&=Cx(t)+Du(t),
\end{aligned}
\end{equation}
where $A,A_j\in\R^{n\times n}$, $B,B_j\in\R^{n\times m}$,
$C\in\R^{m\times n}$, and $D\in\R^{m\times m}$. For a quadratic
storage function
\[
H_X(x)=\tfrac12x^\top Xx,
\qquad X=X^\top,
\]
stochastic passivity with respect to the supply rate $w(u,y)=u^\top y$ is described
by the Kalman--Yakubovich--Popov (KYP) inequality
\begin{equation}
\label{eq:passivity-kyp}
X\succeq0,
\quad
\mathcal M_{\rm pas}(X):=
\begin{bmatrix}
A^\top X+XA+\displaystyle\sum_{j=1}^qA_j^\top XA_j
&
XB-C^\top+\displaystyle\sum_{j=1}^qA_j^\top XB_j
\\[1mm]
B^\top X-C+\displaystyle\sum_{j=1}^qB_j^\top XA_j
&
-(D+D^\top)+\displaystyle\sum_{j=1}^qB_j^\top XB_j
\end{bmatrix}
\preceq0.
\end{equation}
The lower-right block of \eqref{eq:passivity-kyp} contains both the symmetric
feedthrough matrix $D+D^\top$ and the contribution of the
control-dependent diffusion.

A corresponding stochastic available-storage theory for systems with
multiplicative noise has recently been developed in
\cite{AckermannKruseTappe2026}. Under suitable hypotheses, the available storage
is characterized by the Loewner-minimal positive semidefinite solution of
\eqref{eq:passivity-kyp}. Our computational problem is to determine this
minimal solution.

Two features make this problem substantially different from a standard
Riccati computation. First, the symmetric feedthrough matrix $D+D^\top$ may
be singular. In this case Schur complementation cannot be applied directly,
and the singular input directions impose additional structural constraints on
every feasible storage matrix. Second, if the control enters the diffusion,
then the lower-right block of \eqref{eq:passivity-kyp} depends on the unknown
storage matrix. Even after all singular directions have been removed, the
resulting stochastic Riccati equation is therefore a rational matrix equation.

Recursive reduction of singular deterministic KYP problems is developed in
\cite{PoloniReis2012,Reis2011,WeissWangSpeyer1994}. Rational stochastic
Riccati equations and Newton methods for them are treated in
\cite{AitRamiZhou2000,Damm2004,DammHinrichsen2001}. The problem considered
here combines both difficulties: singularity of the constant input block and
control-dependent multiplicative noise. In addition, the relevant solution is
not merely an arbitrary Riccati solution, but the minimal KYP solution
representing the available storage.

We develop a structure-preserving algorithm that keeps this order-theoretic
characterization intact throughout the reduction. We first apply a negative
output-feedback transformation. It leaves the KYP solution set unchanged and,
under the passivity and observability hypotheses used below, yields
exponentially mean-square stable zero-input dynamics. Singular input
directions are then eliminated by exact structural deflation. Feasibility of
the KYP inequality forces linear equations on the storage matrix along these
directions. Solving these equations gives an affine parametrization of all
admissible storage matrices and an induced reduced KYP problem. The reduction
is an order isomorphism, so the minimal solution is transported exactly
through every deflation step.

Stability information is propagated through the same recursion. Starting from
the stabilizing negative output feedback, we transport both a stabilizing
feedback and a strict stochastic Lyapunov certificate to each reduced problem.
Since every nontrivial deflation step decreases the combined state and input
dimension, the process terminates after finitely many steps. The nontrivial
terminal problem then has a positive definite constant input block.

Under a regular-feasibility assumption, the terminal KYP inequality can be
reduced by Schur complementation to a rational stochastic Riccati equation.
Its relevant extremal solution is precisely the deflated available-storage
matrix. A single stochastic Lyapunov solve provides a stabilizing initial point
for Newton's method. The resulting iteration converges monotonically to the
minimal solution; if the limiting feedback is exponentially mean-square
stabilizing, the local convergence is quadratic. Backward reconstruction
through the deflation recursion then recovers the available-storage matrix of
the original system.

For large sparse systems, the reduced Riccati equation is represented directly
on the projected state space. Newton steps lead to generalized stochastic
Lyapunov equations that are represented by signed low-rank factors. Neither a
dense complementary basis nor dense reduced Riccati matrices need to be
assembled. Projected shifted solves retain the sparse state operator, while
the Newton feedback enters only through low-rank corrections. The generalized
Lyapunov equations are solved by a fixed-point iteration whose inner Lyapunov
problems are treated by the signed low-rank alternating-direction implicit (ADI) method. Residual evaluation and backward
reconstruction are performed factorwise as well.

The numerical experiments first compare the dense and low-rank realizations
for a stochastic chain with $n=100$ and then study the scalability of the
sparse implementation for the same problem class. For $n=100$, the dense and
sparse solutions agree to approximately $10^{-8}$. In computations up to
$n=5000$, Newton requires five updates throughout, while the compressed
terminal factors remain of essentially constant width. This illustrates that
the low-rank representation can remain small even though the exact
available-storage matrix is positive definite and hence has full algebraic
rank.

The paper is organized as follows. Section~\ref{sec:preliminaries} recalls the
available-storage result and introduces the general stochastic KYP form used
for the reduction. Section~\ref{sec:deflation} develops the stabilizing loop
transformation, structural deflation, order preservation, and stabilizer
transport. Section~\ref{sec:riccati} treats the terminal stochastic Riccati
equation and Newton iteration, and Section~\ref{sec:algorithm} summarizes the
complete method. Section~\ref{sec:numerics} describes the sparse low-rank
realization, and Section~\ref{sec:experiments} presents the numerical
experiments.

\section{Stochastic KYP inequalities and available storage}
\label{sec:preliminaries}

Throughout this work, we consider a filtered probability space $(\Omega,\mathcal{F},(\mathcal{F}_t)_{t\ge 0},\mathbb{P})$ with a $q$-dimensional Brownian motion $W$ and assume that $(\mathcal{F}_t)_{t\ge 0}$ is the natural filtration of $W$ augmented by all the $\mathbb{P}$-null sets in $\mathcal{F}$. 

Moreover, we introduce the following notation that we will use throughout this work. 

\paragraph{Notation.}
We denote by $\mathbb S^n$, $\mathbb S_+^n$, and $\mathbb S_{++}^n$ the real
symmetric, positive semidefinite, and positive definite $n\times n$ matrices.
Throughout, we allow matrices with zero rows or columns and use the usual
dimensionally consistent conventions for matrix operations. In particular,
the unique $0\times0$ symmetric matrix is regarded as positive definite.
For $X,Y\in\mathbb S^n$, $X\preceq Y$ means $Y-X\succeq0$. The range and
kernel of a matrix $M$ are denoted by $\im M$ and $\ker M$, and $M^\dagger$ denotes
the Moore--Penrose inverse. 
The spectrum and the spectral radius of a bounded linear operator $\mathcal T$ are denoted by $\sigma(\mathcal T)$ and $\rho(\mathcal{T})$. 

\bigskip

Furthermore, we refer to the following 
notions. 

\begin{definition}
The stochastic differential equation 
\begin{equation}\label{eq:uncontrolledsys}
dx(t)=Ax(t)\,dt+\sum_{j=1}^qA_jx(t)\,dW_j(t)
\end{equation}
is called exponentially mean-square stable if 
there exist $c,\gamma >0$ such that 
for all $t\ge 0$ it holds 
$\mathbb E[\|x(t)\|^2] \le ce^{-\gamma t}\|x(0)\|^2$.
\end{definition}

Note that exponential mean-square stability of \eqref{eq:uncontrolledsys} is equivalent to the existence of
$P\in\mathbb S_{++}^n$ with 
\[
A^\top P+PA+\sum_{j=1}^qA_j^\top PA_j \prec 0
\]
(see, e.g., \cite[Theorem~3.6.1(i)]{Damm2004}).



\begin{definition}\label{def:observable}
The stochastic system 
\eqref{eq:system} is called observable 
if for $x$ of \eqref{eq:uncontrolledsys} 
we have that $Cx(t)=0$ $\mathbb{P}$-a.s.\ for all $t\ge 0$ implies $x(0)=0$. 
\end{definition}

Informally, the stochastic system is observable if in the case of zero input and zero output, the state is also zero. 
A convenient sufficient condition for this observability property is the
ordinary deterministic observability condition
\[
\rank\begin{bmatrix}C^\top&(CA)^\top&\cdots&(CA^{n-1})^\top\end{bmatrix}^\top=n;
\]
see, e.g.,  \cite[Remarks~8.8 and~8.11]{AckermannKruseTappe2026}.

Several notions of detectability arise in the stochastic setting;
see \cite{Damm2007}. Here we use the following notion of
$\beta$-detectability.
\begin{definition}\label{def:detectable}
The stochastic system 
\eqref{eq:system} is called $\beta$-detectable 
if for $x$ of \eqref{eq:uncontrolledsys} 
we have that 
$\mathbb{E}[\lVert Cx(t) \rVert^2]=0$ for all $t\ge 0$ 
implies $\lim_{t\to\infty} \mathbb{E}[\lVert x(t) \rVert^2 ]=0$. 
\end{definition}

\subsection{Available storage}


We consider progressively measurable controls that are locally square
integrable. 
More precisely, we denote by $\mathcal{A}$
the set of admissible controls for the stochastic system \eqref{eq:system} 
and introduce the following assumption.

\begin{assumption}\label{assump:admissible_controls}
    The set of admissible controls $\mathcal{A}$ 
    is the set of all progressively measurable processes $u\colon \Omega \times [0,\infty) \to \R^m$ that satisfy 
    \begin{equation*}
        \mathbb{E}\bigg[ \int_0^t \lVert u(s) \rVert^2 \, ds \bigg] < \infty \quad \forall\, t \ge 0 .
    \end{equation*}
\end{assumption}

For the stochastic system \eqref{eq:system}, the available storage
associated with the supply rate $w(u,y)=u^\top y$ is
\begin{equation}
\label{eq:available-storage}
H_{\rm a}(x)=
\sup_{T>0}\sup_{u\in\mathcal{A}}
\mathbb E\!\left[-\int_0^T u(t)^\top y(t)\,dt\right].
\end{equation}
We remark that $H_{\rm a}$ is nonnegative, but not necessarily finite. 
The quadratic storage functions relevant below are precisely those whose
matrices satisfy the passivity KYP inequality \eqref{eq:passivity-kyp}. The
next proposition collects the consequences of
\cite[Proposition~8.4, Remark~8.5, and Lemma~8.13]{AckermannKruseTappe2026}
that are used below. 
In particular, note that \cite[Remark~8.5]{AckermannKruseTappe2026} allows to replace the stochastic passivity assumption in \cite[Proposition~8.4]{AckermannKruseTappe2026} 
by finiteness of \eqref{eq:available-storage}. 


\begin{proposition}
\label{prop:available-storage}
Let  Assumption~\ref{assump:admissible_controls} be satisfied 
and assume that for all $x \in \R^n$ we have $H_{\rm a}(x)<\infty$. 
Then there exists a Loewner-minimal positive semidefinite solution $X_-$ of the stochastic
passivity KYP inequality \eqref{eq:passivity-kyp},  
and it holds for all $x\in\R^n$ that 
\[
H_{\rm a}(x)=\tfrac12x^\top X_-x .
\]
If, in addition, the stochastic system
is observable (in the sense of Definition~\ref{def:observable}), then $X_-\succ0$.
\end{proposition}



\subsection{General stochastic KYP form}\label{sec:generalKYPform}

The deflation analysis is formulated for a general stochastic KYP inequality
with additional constant state and cross terms. Let $Q=Q^\top$,
$S\in\R^{n\times m}$, and $R=R^\top$,
and define
\begin{equation}
\label{eq:general-kyp}
\mathcal M(X)=
\begin{bmatrix}
\mathcal H(X)&\mathcal G(X)\\
\mathcal G(X)^\top&-R_X
\end{bmatrix}\preceq0,
\qquad X\succeq0,
\end{equation}
where
\[
\begin{aligned}
\mathcal H(X)&=A^\top X+XA+\sum_{j=1}^qA_j^\top XA_j+Q,\\
\mathcal G(X)&=XB+\sum_{j=1}^qA_j^\top XB_j+S,\qquad\qquad
R_X=R-\sum_{j=1}^qB_j^\top XB_j.
\end{aligned}
\]
We refer to \eqref{eq:general-kyp} as a stochastic KYP inequality. The
passivity case is recovered from $Q=0$, $S=-C^\top$, and $R=D+D^\top$. The
distinction between the constant input matrix $R$ and the effective matrix
$R_X$ will be crucial below.

\section{Passivity-based loop transformation and structural deflation}
\label{sec:deflation}

\subsection{Stabilizing loop transformation}

We apply a stabilizing loop transformation to the stochastic system \eqref{eq:system} and we also abstract from this to the general KYP problem (see Section~\ref{sec:generalKYPform}). 
To introduce the 
transformation, choose $K=K^\top\succ0$ and assume that $D+D^\top\succeq 0$. 
Note that $D+D^\top\succeq 0$ is guaranteed if the KYP inequality \eqref{eq:passivity-kyp} is feasible. 
It follows that $I+K^{1/2}DK^{1/2}$ 
has positive definite symmetric part and is nonsingular. Since $I+KD$ is
similar to this matrix and it holds $I+DK=K^{-1}(I+KD)K$, both matrices are
nonsingular.
We can thus introduce 
\begin{equation}\label{eq:loop-transform-1}
    L=(I+KD)^{-1},
    \qquad
    F=-LKC 
\end{equation}
and transform the stochastic system \eqref{eq:system} by 
\begin{equation}
\label{eq:loop-transform-2} 
u=Fx+Lv, 
\end{equation}
or, equivalently,  
\[ 
u=v-Ky.
\]
Define for all $j\in\{1,\ldots,q\}$ the matrices 
\begin{align*}
\widehat A&=A+BF,&\widehat B&=BL,&\widehat C & = C + DF,
\\
\widehat A_j&=A_j+B_jF,&
\widehat B_j&=B_jL, & \widehat D &= DL 
\end{align*}
and consider the loop-transformed system 
\begin{equation}
\label{eq:systemlooptrafo}
\begin{aligned}
dx(t)&=(\widehat Ax(t) + \widehat Bv(t))\,dt
 +\sum_{j=1}^q(\widehat A_jx(t)+ \widehat B_jv(t))\,dW_j(t),\\
y(t)&=\widehat C x(t) + \widehat D v(t).
\end{aligned}
\end{equation}
Moreover, note that $\widehat C= (I+DK)^{-1}C$. 
Under the assumptions of Proposition~\ref{prop:available-storage} including the observability assumption, the zero-input dynamics of the loop-transformed system are stable. 

\begin{proposition}
\label{prop:negative-feedback-stable}
Let Assumption~\ref{assump:admissible_controls} be satisfied, suppose that the stochastic system 
\eqref{eq:system} 
is observable, 
and assume that for all $x \in \R^n$ we have $H_{\rm a}(x)<\infty$. 
Then the stochastic differential equation 
\begin{equation}\label{eq:uncontrolledsyshat}
dx(t)=\widehat Ax(t)\,dt+\sum_{j=1}^q \widehat A_jx(t)\,dW_j(t)
\end{equation}
is exponentially mean-square stable.
\end{proposition}

\begin{proof}
We first show that the stochastic system \eqref{eq:systemlooptrafo} is $\beta$-detectable in the sense of Definition~\ref{def:detectable}.
To this end, assume that $x_0\in\R^n$ and that for $x$ of \eqref{eq:uncontrolledsyshat} with initial value $x(0)=x_0$ we have $\mathbb{E}[\lVert \widehat{C} x(t) \rVert^2]=0$ for all $t\ge 0$. 
It follows that $Cx(t)=0$ $\mathbb{P}$-a.s.\ for all $t\ge 0$ and $Fx(t)=0$ $\mathbb{P}$-a.s.\ for all $t\ge 0$.
The latter entails that $x$ solves the original zero-input system
\eqref{eq:uncontrolledsys} with initial value $x_0$,
and using the former and observability of \eqref{eq:system}, we now obtain that $x_0=0$. 
This implies that $x(t)=0$ $\mathbb{P}$-a.s.\ for all $t\ge 0$. We conclude that \eqref{eq:systemlooptrafo} is $\beta$-detectable. 

Let $X$ be a positive definite solution of \eqref{eq:passivity-kyp}, which exists by Proposition~\ref{prop:available-storage}. 
We then have for all $z \in \R^n$ that 
\begin{align}
        \nonumber
        0 & \geq \begin{bmatrix}
            z\\
            Fz
        \end{bmatrix}^\top \mathcal M_{\rm pas}(X) \begin{bmatrix}
            z\\
            Fz
        \end{bmatrix} \\
        & = \nonumber
        z^\top \widehat{A}^\top X z + z^\top X \widehat{A} z + z^\top \sum_{j=1}^q \widehat{A}_j^\top X \widehat{A}_j z
        - z^\top C^\top F z - z^\top F^\top C z
        - z^\top F^\top (D+D^\top) F z \\
        & = 
        z^\top \widehat{A}^\top X z + z^\top X \widehat{A} z + z^\top \sum_{j=1}^q \widehat{A}_j^\top X \widehat{A}_j z 
        + 2 z^\top \widehat{C}^\top K \widehat C z . 
        \label{eq:feedback-lyap-2}
\end{align}
Since $X$ is positive definite and $-2\widehat{C}^\top K \widehat C$ is negative semidefinite, the inequality \eqref{eq:feedback-lyap-2} and $\beta$-detectability of \eqref{eq:systemlooptrafo} show by \cite[Lemma~3.6.6]{Damm2004} that ($\ell$) in \cite[Theorem~3.6.1]{Damm2004} is satisfied. Therefore, \cite[Theorem~3.6.1]{Damm2004} yields that \eqref{eq:uncontrolledsyshat} 
is exponentially mean-square stable.
\end{proof}

Next, we abstract to the general KYP problem. 
In the general KYP case of Section~\ref{sec:generalKYPform}, we assume that $R\succeq 0$, which is guaranteed if the general KYP inequality \eqref{eq:general-kyp} is feasible, and we replace $C$ by $-S^\top$ and $D$ by $\frac12 R$ in \eqref{eq:loop-transform-1}, respectively. 
Furthermore, let 
\[
\widehat Q=Q+SF+F^\top S^\top-F^\top RF, \qquad
\widehat S=(S-F^\top R)L, \qquad \widehat R=L^\top RL,\qquad
T=\begin{bmatrix}I&0\\F&L\end{bmatrix}.
\]
The transformed KYP problem is 
\[
\widehat{\mathcal M}(X) \preceq 0, \quad X\succeq 0,
\]
where $\widehat{\mathcal M}(X)$ is defined analogously to $\mathcal M(X)$, but with the coefficient matrices $A$, $A_j$, $B$, $B_j$, $Q$, $S$, $R$ replaced by the transformed KYP data $\widehat A, \widehat A_j, \widehat B, \widehat B_j, \widehat Q, \widehat S, \widehat R$. 
The transformation procedure preserves the solution set of the KYP inequality and the rank of $R$. 

\begin{proposition}\label{prop:transformed-KYP-equivalent}
Assume that $R\succeq 0$. 
Let $F$, $L$, $T$, and the transformed KYP data be defined as above. Then it holds for all $X=X^\top$ that 
\[
\widehat{\mathcal M}(X)=T^\top\mathcal M(X)T.
\]
Consequently, 
\[
\{X \,|\, \mathcal{M}(X) \preceq 0, \, X\succeq 0 \}
= \{X \,|\, \widehat{\mathcal{M}}(X) \preceq 0, \, X\succeq 0 \}, 
\]
and
$\rank\widehat R=\rank R$.
\end{proposition}

\begin{proof}
The first assertion follows by direct multiplication. Since $T$ is invertible,
negative semidefiniteness is equivalent in both representations. The rank
identity follows from the invertibility of $L$.
\end{proof}

In the deflation analysis we suppress the hats and write
$A,A_j,B,B_j,Q,S,R$ for the current loop-transformed KYP data.

\subsection{One deflation step}\label{sec:onedeflationstep}

Let $R\in\mathbb S_+^m$ be singular and set $s=\dim\ker R$.
Choose an orthogonal matrix $U=[U_1\ U_0]$, where
$U_1\in\R^{m\times(m-s)}$ and $U_0\in\R^{m\times s}$,
such that
\[
U^\top RU=
\begin{bmatrix}R_1&0\\0&0\end{bmatrix},
\qquad R_1\in\mathbb S_{++}^{m-s}.
\]
For all $j\in\{1,\ldots,q\}$ partition
\[
BU=[B_{\rm r}\ B_0],\qquad
B_jU=[B_{j1}\ B_{j0}],\qquad
SU=[S_1\ S_0].
\]
Here $B_{\rm r},B_{j1},S_1\in\R^{n\times(m-s)}$ and
$B_0,B_{j0},S_0\in\R^{n\times s}$.
Set
\begin{equation}\label{eq:defEJ}
E=[B_0\ B_{10}\ \cdots\ B_{q0}],
\qquad
J=[-S_0\ 0\ \cdots\ 0].
\end{equation}
Both $E$ and $J$ belong to $\R^{n\times(q+1)s}$.
The next lemma implies that every feasible matrix $X$ for the stochastic KYP problem \eqref{eq:general-kyp} satisfies $XE=J$.

\begin{lemma}
\label{lem:singular-directions}
Consider a feasible stochastic KYP problem \eqref{eq:general-kyp} with
singular $R\succeq0$. 
Then every $X\succeq0$ satisfying \eqref{eq:general-kyp} fulfills
\[
XB_0=-S_0,
\qquad
XB_{j0}=0,
\quad j=1,\ldots,q.
\]
\end{lemma}

\begin{proof}
Suppose that $X\succeq0$ satisfies \eqref{eq:general-kyp} and let $z\in\ker R$. 
The lower-right block of the negative semidefinite matrix $\mathcal{M}(X)$ then gives
\[
0\le z^\top R_Xz
=-\sum_{j=1}^q(B_jz)^\top X(B_jz)\le0,
\]
hence $XB_jz=0$ for all $j\in\{1,\ldots,q\}$. Since the quadratic form of the negative semidefinite matrix
$\mathcal M(X)$ vanishes on $[0\ z^\top]^\top$, this vector belongs to
$\ker\mathcal M(X)$. Thus $\mathcal G(X)z=0$ and hence $XBz+Sz=0$.
\end{proof}


The next proposition characterizes 
the set of solutions to the equations $X=X^\top\succeq0$ and $XE=J$. 

\begin{proposition}
\label{prop:affine-parametrization}
Let $R\succeq 0$ be singular 
and let $E$ and $J$ be defined as in \eqref{eq:defEJ}. 
Let $V$ be a matrix whose columns form an orthonormal basis of $\ker E^\top$. 
Then the equations 
\begin{equation}\label{eq:eqXEJ}
    X=X^\top\succeq0, \qquad XE=J
\end{equation}
are solvable if and only if it holds 
\begin{equation}\label{eq:consistencycondEJ}
E^\top J=(E^\top J)^\top\succeq0,
\qquad
\ker(E^\top J)=\ker J.
\end{equation}
In this case
\[
X_{\rm p}=J(E^\top J)^\dagger J^\top
\]
is the Loewner-minimal positive semidefinite solution 
of \eqref{eq:eqXEJ}, and it holds  
\begin{equation*}
    \{X \,|\, X\succeq 0, \, XE=J \}
    = \{X \,|\, \exists Z\succeq 0\colon X=X_{\rm p}+VZV^\top\} .
\end{equation*}
\end{proposition}

\begin{proof}
Necessity 
of \eqref{eq:consistencycondEJ} 
follows from $E^\top J=E^\top XE\succeq0$ and
$\ker(E^\top XE)=\ker(XE)=\ker J$. 
Conversely, assume \eqref{eq:consistencycondEJ} and set $H=E^\top J$.
Since $H\succeq0$ and $\ker H=\ker J$, one has
$J=JH^\dagger H$. Hence
\[
X_{\rm p}=JH^\dagger J^\top\succeq0,
\qquad
X_{\rm p}E=JH^\dagger H=J,
\]
i.e., $X_{\rm p}$ is a solution of \eqref{eq:eqXEJ}. 
If $X\succeq0$ and $XE=J$, then
\[
X_{\rm p}=XE(E^\top XE)^\dagger E^\top X\preceq X,
\]
because, after writing $X=X^{1/2}X^{1/2}$, the middle factor is the
orthogonal projector onto $\im(X^{1/2}E)$. 
Moreover,
$(X-X_{\rm p})E=0$, so the symmetric positive semidefinite matrix
$X-X_{\rm p}$ has range in $\ker E^\top=\im V$.  
Therefore, $Z:=V^\top (X-X_{\rm p}) V$ is symmetric positive semidefinite and satisfies $X-X_{\rm p}=VZV^\top$. 
\end{proof}

In the sequel, let $r=\rank E$ and fix matrices
$V\in\R^{n\times(n-r)}$ and $\bar V\in\R^{n\times r}$
whose columns form orthonormal bases of $\ker E^\top$
and $\im E$, respectively.
Then $\mathcal V=[V\ \bar V]\in\R^{n\times n}$ is orthogonal. 
Set $X_{\rm p}=J(E^\top J)^\dagger J^\top$ and define 
\begin{align*}
\Theta&=A^\top X_{\rm p}+X_{\rm p}A
+\sum_{j=1}^qA_j^\top X_{\rm p}A_j+Q,\\
\Xi&=X_{\rm p}B_{\rm r}+\sum_{j=1}^qA_j^\top X_{\rm p}B_{j1}+S_1,\qquad\qquad
\Upsilon=-R_1+\sum_{j=1}^qB_{j1}^\top X_{\rm p}B_{j1}.
\end{align*}
Here $\Theta\in\mathbb S^n$,
$\Xi\in\R^{n\times(m-s)}$, and
$\Upsilon\in\mathbb S^{m-s}$.
The reduced data are 
\begin{equation}
\label{eq:reduced-data}
\begin{aligned}
A^+&=V^\top AV,
& A_j^+&=V^\top A_jV,\\
B^+&=[V^\top A\bar V \ \ V^\top B_{\rm r}],
& B_j^+&=[V^\top A_j\bar V\ \ V^\top B_{j1}],\\
Q^+&=V^\top\Theta V,
& S^+&=[V^\top\Theta \bar V\ \ V^\top\Xi],&
R^+&=-\begin{bmatrix}
\bar{V}^\top\Theta \bar V&\bar{V}^\top\Xi\\
\Xi^\top \bar V&\Upsilon
\end{bmatrix}
\end{aligned}
\end{equation}
for all $j\in\{1,\ldots,q\}$.
For $Z\in\mathbb S_+^{n-r}$, the reduced KYP problem is
\begin{equation}\label{eq:KYP-reduced}
\mathcal M^+(Z) \preceq 0, \quad Z\succeq 0,
\end{equation}
where $\mathcal M^+(Z)$ is defined analogously to $\mathcal M(X)$ in \eqref{eq:general-kyp}, but with the coefficient matrices $A$, $A_j$, $B$, $B_j$, $Q$, $S$, $R$ replaced by the reduced KYP data $A^+,A_j^+,B^+,B_j^+,Q^+,S^+,R^+$. 
The reduced state and input dimensions are
$n^+=n-r$ and $m^+=m-s+r$, respectively.
In particular,
$A^+,A_j^+,Q^+\in\R^{n^+\times n^+}$,
$B^+,B_j^+,S^+\in\R^{n^+\times m^+}$,
and $R^+\in\mathbb S^{m^+}$.

\begin{theorem}
\label{thm:deflation-step}
Consider a feasible stochastic KYP problem \eqref{eq:general-kyp} with singular 
$R\succeq0$ 
and define the reduced data
by \eqref{eq:reduced-data}. Then the map
\[
Z\longmapsto X_{\rm p}+VZV^\top
\]
is a bijection between the positive semidefinite feasible sets of the reduced
and original KYP problems. The reduced constant input block satisfies
$R^+\succeq0$. 
\end{theorem}

\begin{proof}
First, observe that feasibility of the stochastic KYP problem \eqref{eq:general-kyp} and Lemma~\ref{lem:singular-directions} ensure that the consistency conditions \eqref{eq:consistencycondEJ} hold. 
Second, note that for any symmetric matrix $Z$ and $X:=X_{\rm p} + V Z V^\top$ one can show using the properties $\mathcal{V}^\top\mathcal{V}=I$, $X_{\rm p} E = J$ (cf.\ Proposition~\ref{prop:affine-parametrization}), and $E^\top V = 0$ that 
\begin{equation}\label{eq:KYPmatEquiv}
\begin{bmatrix}\mathcal V&0\\0&U\end{bmatrix}^\top
\mathcal M(X)
\begin{bmatrix}\mathcal V&0\\0&U\end{bmatrix}
\sim
\begin{bmatrix}
\mathcal M^+(Z)&0\\
0&0
\end{bmatrix}
\end{equation}
(where $\sim$ means up to some permutation).
If $Z$ is a solution of the reduced KYP problem~\eqref{eq:KYP-reduced}, then $X_{\rm p} + V Z V^\top$ is a solution of the original KYP problem~\eqref{eq:general-kyp} by Proposition~\ref{prop:affine-parametrization} and \eqref{eq:KYPmatEquiv}. 
Next, suppose that $X$ solves the original KYP problem~\eqref{eq:general-kyp}. Then Lemma~\ref{lem:singular-directions} guarantees that $X$ satisfies \eqref{eq:eqXEJ}, and hence by Proposition~\ref{prop:affine-parametrization} there exists $Z\succeq 0$ such that $X=X_{\rm p} + V Z V^\top$. Since $\mathcal{M}(X)\preceq 0$, it follows from \eqref{eq:KYPmatEquiv} that $Z$ satisfies~\eqref{eq:KYP-reduced}.   
Moreover, note that if $Z_1,Z_2$ satisfy~\eqref{eq:KYP-reduced} and $X_{\rm p} + V Z_1 V^\top=X_{\rm p} + V Z_2 V^\top$, then $Z_1=Z_2$ due to $V^\top V = I$.
Furthermore, for any solution $Z$ of~\eqref{eq:KYP-reduced}, the lower-right block $-R_Z^+$ of $\mathcal{M}^+(Z)$ gives
\[
0 \preceq R_Z^+=R^+-\sum_{j=1}^q(B_j^+)^\top ZB_j^+. 
\]
Since the original problem is feasible, so is the reduced one, and hence $R^+\succeq0$.
\end{proof}



\begin{remark}[Degenerate dimensions]
The preceding construction also covers degenerate dimensions
without requiring separate deflation formulas.
If $R=0$, the block $R_1$ is empty, while the constraints
$XE=J$ remain in force.
If $E=0$, feasibility implies $J=0$, and the singular input
directions are eliminated without reducing the state dimension.
If $\rank E=n$, the reduced state dimension is zero, and
the recursion terminates.
Finally, if $R\succ0$, no deflation is required.
Thus all cases are covered by the same construction and the
terminal conditions of the algorithm.
\end{remark}

\subsection{Finite deflation, order, and factor reconstruction}\label{sec:finitedeflation}

At deflation level $k$ we write all data with a superscript $(k)$. Thus,
$R^{(k)}$ denotes the constant input matrix of the current KYP problem, whereas
the subscript $0$ denotes its singular input directions (as in Section~\ref{sec:onedeflationstep}). If
$s_k=\dim\ker R^{(k)}$ and $r_k=\rank E^{(k)}$, then
\[
n_{k+1}+m_{k+1}=n_k+m_k-s_k.
\]
Every nontrivial deflation step decreases $n_k+m_k$.
Hence the procedure terminates after finitely many steps, stopping
when $n_\ell=0$, $m_\ell=0$, or $R^{(\ell)}\succ0$.

The parametrization is an order isomorphism:
\begin{equation}
\label{eq:order-equivalence}
Z_1\preceq Z_2
\quad\Longleftrightarrow\quad
X_{\rm p}+VZ_1V^\top\preceq X_{\rm p}+VZ_2V^\top.
\end{equation}
Consequently, if $X_-^{(k)}$ denotes the minimal KYP solution at level $k$,
then
\begin{equation}
\label{eq:minimal-reconstruction}
X_-^{(k)}=X_{{\rm p}}^{(k)}
+V^{(k)}X_-^{(k+1)}(V^{(k)})^\top.
\end{equation}
If $X_-^{(k)}\succ0$, then $X_-^{(k+1)}\succ0$: indeed,
$\rank X_{\rm p}^{(k)}\le r_k$, whereas singularity of
$X_-^{(k+1)}$ would force the right-hand side of
\eqref{eq:minimal-reconstruction} to have rank at most $n_k-1$.

Backward reconstruction is represented factorwise. Let
\[
H_k=(E^{(k)})^\top J^{(k)},
\qquad
Z_{{\rm p},k}=J^{(k)}(H_k^\dagger)^{1/2}.
\]
Then $X_{{\rm p}}^{(k)}=Z_{{\rm p},k}Z_{{\rm p},k}^\top$. If
$X_-^{(k+1)}=Z_{k+1}Z_{k+1}^\top$, then
\begin{equation}
\label{eq:factor-reconstruction}
X_-^{(k)}=Z_kZ_k^\top,
\qquad
Z_k=\begin{bmatrix}Z_{{\rm p},k}&V^{(k)}Z_{k+1}\end{bmatrix}.
\end{equation}
The columns of the reconstructed factor $Z_k$ may be recompressed after each
step.

The backward recursion starts from the terminal storage matrix $X_-^{(\ell)}$. The computation of the minimal solution $X_-^{(\ell)}$ of the terminal KYP in the case $R^{(\ell)} \succ 0$ is the subject of Section~\ref{sec:riccati}, while the cases $n_{\ell}=0$ and $m_{\ell}=0$ are discussed in the next two remarks. 

\begin{remark}[Termination without remaining states]
\label{rem:zero-terminal-state}
If $n_\ell=0$, the terminal KYP inequality reduces to
$-R^{(\ell)}\preceq0$, which holds by feasibility.
The terminal storage matrix is the unique $0\times0$ matrix.
Thus no terminal matrix equation needs to be solved, and
backward factor reconstruction starts with an empty factor.
\end{remark}

\begin{remark}[Termination without remaining inputs]
\label{rem:zero-terminal-input}
If $m_\ell=0$, the terminal KYP inequality reduces to
\[
X\succeq0,\qquad
\mathcal L_\ell(X)+Q^{(\ell)}\preceq0,
\]
where
\[
\mathcal L_\ell(X)
=(A^{(\ell)})^\top X+XA^{(\ell)}
+\sum_{j=1}^q(A_j^{(\ell)})^\top XA_j^{(\ell)}.
\]
Suppose that there is a minimal solution $X_-^{(\ell)}\succ0$ of the terminal KYP inequality and that 
the terminal dynamics are exponentially mean-square
stable. 
Then the minimal terminal KYP solution
is the unique solution of
\[
\mathcal L_\ell(X)+Q^{(\ell)}=0.
\]
Indeed, if 
the inequality were not an equality at $X_-^{(\ell)}$ in the sense that
\[
Y_\ell=-\mathcal L_\ell(X_-^{(\ell)})-Q^{(\ell)}
\succeq0,\qquad Y_\ell\ne0,
\]
the stochastic Lyapunov equation
$\mathcal L_\ell(Z)=-Y_\ell$ would have a nonzero
positive semidefinite solution $Z$.
For sufficiently small $\varepsilon>0$,
$X_-^{(\ell)}-\varepsilon Z$ would remain
positive definite and KYP-feasible, contradicting
minimality. Thus $Y_\ell=0$.
Uniqueness follows from exponential mean-square stability.
Consequently, no Riccati iteration is needed in this case.
\end{remark}

\subsection{Transport of a stabilizing feedback}
\label{sec:stabilizertransport}

A feedback matrix $\Phi$ specifies the static state feedback $u=\Phi x$
for the current system.
After the initial loop transformation, the zero feedback is stabilizing (see Proposition~\ref{prop:negative-feedback-stable}). 
Repeated application of Theorem~\ref{thm:stabilizer-transport} below therefore
provides a stabilizing terminal feedback $\Phi^{(\ell)}$.
Note that Theorem~\ref{thm:stabilizer-transport} employs the notation introduced in Section~\ref{sec:onedeflationstep}. 

\begin{theorem}
\label{thm:stabilizer-transport}
Let $R\succeq 0$ be singular. 
Suppose that a feedback $\Phi$ makes
$A_\Phi=A+B\Phi$ and $A_{j,\Phi}=A_j+B_j\Phi$, $j\in\{1,\ldots,q\}$, exponentially mean-square
stable, and let $P\succ0$ satisfy
\[
\mathcal L_\Phi(P):=
A_\Phi^\top P+PA_\Phi+
\sum_{j=1}^qA_{j,\Phi}^\top PA_{j,\Phi}\prec0.
\]
With $\mathcal V=[V\ \bar V]$, partition
\[
\mathcal V^\top P\mathcal V=
\begin{bmatrix}P_{11}&P_{12}\\P_{12}^\top&P_{22}\end{bmatrix}
\]
and set
\[
\Lambda=-P_{22}^{-1}P_{12}^\top,
\qquad
T_P=\begin{bmatrix}I\\\Lambda\end{bmatrix},
\qquad
P^+=P_{11}-P_{12}P_{22}^{-1}P_{12}^\top.
\]
If $\Phi_1=U_1^\top\Phi$, define
\[
\Phi^+=
\begin{bmatrix}
\Lambda\\[1mm]
\Phi_1\mathcal VT_P
\end{bmatrix}.
\]
Then $\Phi^+$ exponentially mean-square stabilizes the reduced stochastic
system, and $P^+\succ0$ is a strict Lyapunov certificate.
\end{theorem}

\begin{proof}
Since $V^\top B_0=V^\top B_{j0}=0$, it holds 
\[
A_{\Phi^+}^+=V^\top A_\Phi\mathcal VT_P,
\qquad
A_{j,\Phi^+}^+=V^\top A_{j,\Phi}\mathcal VT_P.
\]
Moreover,
\[
\mathcal V^\top P\mathcal VT_P
=\begin{bmatrix}P^+\\0\end{bmatrix}.
\]
The Schur-complement identity gives
\[
P-VP^+V^\top\succeq0.
\]
Consequently, the drift terms in the reduced Lyapunov expression agree with
the compression of the corresponding full-order terms, while for every noise
matrix
\[
(A_{j,\Phi^+}^+)^\top P^+A_{j,\Phi^+}^+
\preceq
T_P^\top\mathcal V^\top A_{j,\Phi}^\top P
A_{j,\Phi}\mathcal VT_P.
\]
Summing these relations yields
\[
\mathcal L_{\Phi^+}^+(P^+)
\preceq T_P^\top\mathcal V^\top\mathcal L_\Phi(P)\mathcal VT_P\prec0.
\]
Moreover, $P^+=T_P^\top\mathcal V^\top P \mathcal VT_P \succ 0$. 
\end{proof}


\section{Terminal stochastic Riccati equation and Newton iteration}
\label{sec:riccati}

We suppress the terminal index $\ell$ in this section and consider the
nontrivial terminal case $n_\ell>0$ and $m_\ell>0$. Thus $R\succ0$ and $X_-$ denotes the
minimal positive semidefinite solution of the terminal KYP problem. If the
original passive system is observable, Proposition~\ref{prop:available-storage}
and the positive-definiteness propagation in Section~\ref{sec:deflation} give
$X_-\succ0$.

\subsection{Riccati reduction}

Define the regular domain
\[
\mathcal D=\{X=X^\top\mid R_X\succ0\}.
\]

\begin{assumption}[Regular feasibility]
\label{ass:regular-feasibility}
There exists $\widehat X=\widehat X^\top$ such that
\[
\mathcal M(\widehat X)\preceq0,
\qquad
R_{\widehat X}\succ0.
\]
\end{assumption}

Regular feasibility is an existence hypothesis: the witness $\widehat X$ is
required only in the analysis. For $X\in\mathcal D$ define
\begin{equation}
\label{eq:riccati-map}
\mathcal R(X)=
\mathcal H(X)+\mathcal G(X)R_X^{-1}\mathcal G(X)^\top.
\end{equation}
Then Schur complementation gives
\[
\mathcal M(X)\preceq0
\quad\Longleftrightarrow\quad
\mathcal R(X)\preceq0,
\qquad X\in\mathcal D,
\]
and the terminal stochastic algebraic Riccati equation is
\begin{equation}
\label{eq:stare}
\mathcal R(X)=0.
\end{equation}
For control-dependent noise this is a rational matrix equation of the type
studied in \cite{Damm2004,DammHinrichsen2001}.

For $X\in\mathcal D$ set
\[
\Gamma(X)=R_X^{-1}\mathcal G(X)^\top.
\]
For a feedback $\Phi \in \R^{m\times n}$ and all $j\in\{1,\ldots,q\}$ define
\[
A_\Phi=A+B\Phi,
\qquad
A_{j,\Phi}=A_j+B_j\Phi,
\qquad
Q_\Phi=Q+S\Phi+\Phi^\top S^\top-\Phi^\top R\Phi,
\]
and let $\mathcal L_\Phi$ be the associated stochastic Lyapunov operator defined by 
\[
\mathcal L_\Phi(X):=
A_\Phi^\top X+XA_\Phi+
\sum_{j=1}^qA_{j,\Phi}^\top X A_{j,\Phi}.
\]

\begin{lemma}
\label{lem:riccati-identities}
For $X\in\mathcal D$ and arbitrary $\Phi \in \R^{m\times n}$ it holds 
\begin{equation}
\label{eq:completion-squares}
\mathcal R(X)=
\mathcal L_\Phi(X)+Q_\Phi
+(\Gamma(X)-\Phi)^\top R_X(\Gamma(X)-\Phi),
\end{equation}
and, for every symmetric $Z$, 
\[
D\mathcal R(X)[Z]=\mathcal L_{\Gamma(X)}(Z).
\]
\end{lemma}

\begin{proof}
The first identity is obtained by completing the square in
$\mathcal G(X)R_X^{-1}\mathcal G(X)^\top$. Differentiating
\eqref{eq:riccati-map} and using
$D(R_X^{-1})[Z]=R_X^{-1}(\sum_jB_j^\top ZB_j)R_X^{-1}$ gives the second
identity after collecting the closed-loop terms.
\end{proof}

\subsection{Stabilizing initialization}

The deflation recursion provides an exponentially mean-square
stabilizing feedback for the terminal KYP problem (see Section~\ref{sec:stabilizertransport}). To initialize
Newton's method and to establish stability of the subsequent
Riccati feedbacks, we use the following elementary
stability-transfer result.
\begin{lemma}\label{lem:feedback-stability-transfer}
Let $\Phi$ be an exponentially mean-square stabilizing
feedback, and let $\Psi$ be another feedback.
Suppose that there exist $Y\succeq0$ and $M\succ0$
such that
\[
\mathcal L_\Psi(Y)\preceq
-(\Psi-\Phi)^\top M(\Psi-\Phi).
\]
Then $\Psi$ is exponentially mean-square stabilizing.
\end{lemma}

\begin{proof}
Choose $P\succ0$ with $\mathcal L_\Phi(P)\prec0$
and set $N=M^{1/2}(\Psi-\Phi)$.
For $z\in\ker N$, the closed-loop drift and diffusion
matrices associated with $\Psi$ and $\Phi$ agree
when applied to $z$. Consequently,
\[
z^\top\mathcal L_\Psi(P)z
=z^\top\mathcal L_\Phi(P)z<0
\]
for every nonzero $z\in\ker N$.
It follows that there exists $\tau>0$ such that
\[
\mathcal L_\Psi(P)-\tau N^\top N\prec0.
\]
Indeed, otherwise a sequence of unit vectors corresponding
to $\tau\to\infty$ would have a subsequence converging
to a nonzero vector in $\ker N$ on which
$\mathcal L_\Psi(P)$ is nonnegative, a contradiction.
Hence
\[
\mathcal L_\Psi(P+\tau Y)
\preceq\mathcal L_\Psi(P)-\tau N^\top N\prec0.
\]
Since $P+\tau Y\succ0$, the feedback $\Psi$
is exponentially mean-square stabilizing.
\end{proof}

For an exponentially mean-square stabilizing terminal feedback $\Phi$, the
stochastic Lyapunov equation is
\begin{equation}
\label{eq:initial-feedback-lyapunov}
\mathcal L_\Phi(X_0)+Q_\Phi=0.
\end{equation}
By stability, there exists a unique symmetric
solution $X_0$ of \eqref{eq:initial-feedback-lyapunov}. 
The next result shows, under
Assumption~\ref{ass:regular-feasibility}, that $X_0\in\mathcal D$ and that $X_0$ can serve as a stabilizing initial matrix for \eqref{eq:stare}. 
The initial matrix $X_0$ need not be positive semidefinite.

\begin{theorem}
\label{thm:X0-stabilizing}
Let $\Phi\in\R^{m\times n}$ be any exponentially mean-square stabilizing
feedback, and let $X_0$ be the unique symmetric solution of
\eqref{eq:initial-feedback-lyapunov}.
Under Assumption~\ref{ass:regular-feasibility}, it holds that
\[
X_0\preceq\widehat X,\qquad
R_{X_0}\succeq R_{\widehat X}\succ0,\qquad
\mathcal R(X_0)\succeq0.
\]
Moreover, the Riccati feedback $\Gamma_0:=\Gamma(X_0)$ is exponentially
mean-square stabilizing.
\end{theorem}


\begin{proof}
Set $Y=\widehat X-X_0$. Evaluating the quadratic form of the
KYP matrix $\mathcal M(\widehat X)\preceq0$ along the feedback $\Phi$
gives $\mathcal L_\Phi(\widehat X)+Q_\Phi\preceq0$.
Subtracting \eqref{eq:initial-feedback-lyapunov} yields
$\mathcal L_\Phi(Y)\preceq0$. Since $\Phi$ is exponentially mean-square
stabilizing, the integral representation of the stochastic Lyapunov equation yields
$Y\succeq0$. Thus, 
\[
R_{X_0}=R_{\widehat X}+\sum_{j=1}^qB_j^\top YB_j\succeq R_{\widehat X}.
\]
Hence, using Assumption~\ref{ass:regular-feasibility},  $X_0\in\mathcal D$. 
With $\Delta=\Gamma_0-\Phi$, \eqref{eq:completion-squares} gives
$\mathcal R(X_0)=\Delta^\top R_{X_0}\Delta\succeq0$. Applying the same identity
at $\widehat X$ with the fixed feedback $\Gamma_0$ and recalling that, due to Assumption~\ref{ass:regular-feasibility}, it holds $R_{\widehat X} \succ 0$ and $\mathcal R(\widehat X)\preceq 0$ yields
\[
\mathcal L_{\Gamma_0}(Y)
=\mathcal R(\widehat X)-\mathcal R(X_0)-
\bigl(\Gamma(\widehat X)-\Gamma_0\bigr)^\top
R_{\widehat X}
\bigl(\Gamma(\widehat X)-\Gamma_0\bigr)\preceq-\Delta^\top R_{X_0}\Delta.
\]
Since $Y\succeq0$ and $R_{X_0}\succ0$,
Lemma~\ref{lem:feedback-stability-transfer}, applied with
$\Psi=\Gamma_0$ and $M=R_{X_0}$, shows that
$\Gamma_0$ is exponentially mean-square stabilizing.
\end{proof}





\subsection{Extremal solution and Newton convergence}

The next result identifies the minimal terminal KYP solution with the
least solution of the regular Riccati inequality. We prove convergence of
Newton's method directly using regular feasibility and the stabilizing
initial feedback, without an additional change of sign convention.

By Lemma~\ref{lem:riccati-identities}, Newton's equation 
\[
D\mathcal R(X_k) [X_{k+1}-X_k] = - \mathcal R(X_k)
\]
is equivalent to the
stochastic Lyapunov equation
\begin{equation}
\label{eq:newton-lyapunov}
\mathcal L_{\Gamma_k}(X_{k+1})+Q_{\Gamma_k}=0,
\qquad
\Gamma_k=R_{X_k}^{-1}\mathcal G(X_k)^\top.
\end{equation}
Thus every Newton step is a stochastic Lyapunov solve.

\begin{theorem}
\label{thm:riccati-storage-identification}
Assume $X_-\succ0$ for the minimal solution of $\mathcal{M}(X)\preceq 0$, $X\succeq 0$, and let Assumption~\ref{ass:regular-feasibility} be in force. 
Let $\Phi \in \R^{m\times n}$ be exponentially mean-square stabilizing 
and let $X_0$ be the unique symmetric solution of \eqref{eq:initial-feedback-lyapunov}. 
Then $R_{X_-}\succ0$ and $\mathcal R(X_-)=0$. Moreover, $X_-$ is the least
solution of the regular Riccati inequality $\mathcal R(X)\preceq0$
on $\mathcal D$, without any positivity restriction on $X$.
Newton's method initialized with $X_0$ 
is well defined and converges monotonically to $X_-$:
\[
X_0\preceq X_1\preceq X_2\preceq\cdots\preceq X_-,
\qquad X_k\longrightarrow X_-.
\]
Each Riccati feedback $\Gamma_k:=\Gamma(X_k)$ is exponentially mean-square stabilizing.
At the limit, $D\mathcal R(X_-)$ has spectrum in the closed left half-plane.
\end{theorem}

\begin{proof}
Since $X_{k+1}$ solves \eqref{eq:initial-feedback-lyapunov}
with $\Phi=\Gamma_k$, Theorem~\ref{thm:X0-stabilizing}
applies inductively. Thus all iterates belong to $\mathcal D$,
have stabilizing Riccati feedbacks, and satisfy
$X_k\preceq\widehat X$,
$R_{X_k}\succeq R_{\widehat X}\succ0$, and
$\mathcal R(X_k)\succeq0$ for all $k\ge0$. 
Note that 
\[
\mathcal L_{\Gamma_k}(X_{k+1}-X_k)
=-\mathcal R(X_k)\preceq0
\]
and exponential mean-square stability 
imply $X_{k+1}\succeq X_k$. Evaluating the quadratic form of $\mathcal M(X_-)\preceq0$
along $\Gamma_k$ gives
$\mathcal L_{\Gamma_k}(X_-)+Q_{\Gamma_k}\preceq0$.
Comparison with \eqref{eq:newton-lyapunov} and exponential
mean-square stability yields $X_{k+1}\preceq X_-$. Thus $(X_k)$ converges in the Loewner order
to some $X_\infty\preceq X_-$. Since
$X_k\preceq X_{k+1} \preceq \widehat X$ for all $k$ and
$R_{X_k}\succeq R_{X_{k+1}}\succeq R_{\widehat X}\succ0$, the feedbacks
$\Gamma_k$ and Lyapunov operators $\mathcal L_{\Gamma_k}$
are uniformly bounded and it holds $X_{\infty} \in \mathcal D$. As $X_{k+1}-X_k\to0$, 
the equation $\mathcal L_{\Gamma_k}(X_{k+1}-X_k)
=-\mathcal R(X_k)$ 
implies $\mathcal R(X_k)\to0$ 
and thus $\mathcal R(X_\infty)=0$. 
By Schur complementation, this and $X_\infty\in\mathcal D$ yield $\mathcal M(X_\infty)\preceq0$.

If $X_\infty\ne X_-$, then for sufficiently small $\theta>0$,
\[
X_\theta=(1-\theta)X_-+\theta X_\infty
\]
is positive definite, satisfies $\mathcal M(X_\theta)\preceq0$
by affinity, and obeys $X_\theta\preceq X_-$ with
$X_\theta\ne X_-$. This contradicts the Loewner minimality of
$X_-$. Hence $X_\infty=X_-$. Finally, for any
$\widetilde X\in\mathcal D$ with $\mathcal R(\widetilde X)\preceq0$,
\eqref{eq:completion-squares} yields
$\mathcal L_{\Gamma_k}(\widetilde X)+Q_{\Gamma_k}\preceq0$.
Comparison with the equation for $X_{k+1}$ and exponential
mean-square stability give $X_{k+1}\preceq\widetilde X$.
Letting $k\to\infty$ proves $X_-\preceq\widetilde X$.
Stability of every $\mathcal L_{\Gamma_k}$ and continuity
of the spectrum imply that $\mathcal L_{\Gamma(X_-)}$, and thus $D\mathcal R(X_-)$ by Lemma~\ref{lem:riccati-identities}, 
has spectrum in the closed left half-plane.
\end{proof}

\begin{remark}[Local quadratic convergence]
If, in addition, the limiting feedback
\[
\Gamma_-:=\Gamma(X_-)
\]
is exponentially mean-square stabilizing, then
\[
D\mathcal R(X_-)=\mathcal L_{\Gamma_-}
\]
is invertible. Since $\mathcal R$ is smooth on the regular domain
$\mathcal D$ and $D\mathcal R(X_-)$ is invertible,
the inverses of $D\mathcal R(X_k)$ remain uniformly bounded
for all sufficiently large $k$. Taylor's formula applied
to the Newton equation therefore yields
\[
\|X_{k+1}-X_-\|\le c\|X_k-X_-\|^2
\]
for some $c>0$ and all sufficiently large $k$. This additional stability
assumption is not needed for the monotone convergence result in
Theorem~\ref{thm:riccati-storage-identification}.
\end{remark}


\section{Deflation--Riccati algorithm}
\label{sec:algorithm}

This section summarizes the procedure to compute the available storage of the stochastic system~\eqref{eq:system}. 
As in Section~\ref{sec:finitedeflation}, superscripts $(k)$ index deflation levels. 
The complete procedure is:

\begin{enumerate}[label=\textbf{\arabic*.},leftmargin=2.5em]
\item Choose $K=K^\top\succ0$, apply the loop transformation
      \eqref{eq:loop-transform-1}--\eqref{eq:loop-transform-2} to the stochastic system~\eqref{eq:system},
      and form the equivalent KYP data. Set
      $\Phi^{(0)}=0$ and compute a strict stochastic Lyapunov certificate,
      whose existence follows from Proposition~\ref{prop:negative-feedback-stable}.

\item While $n_k>0$, $m_k>0$, and $R^{(k)}$ is singular,
      perform the deflation of
      Theorem~\ref{thm:deflation-step}, store $Z_{{\rm p},k}$ together with
      a representation of the complementary state subspace (for instance
      $V^{(k)}$, or $\bar V^{(k)}$ and the associated projector), and transport the
      stabilizing feedback and certificate by
      Theorem~\ref{thm:stabilizer-transport}.

\item At level $\ell$, distinguish three cases:
        \begin{itemize}
            \item If $n_\ell=0$, set $X_-^{(\ell)}=0_{0\times0}$ as in Remark~\ref{rem:zero-terminal-state}.
            
            \item If $n_\ell>0$ and $m_\ell=0$, solve the stochastic Lyapunov equation $\mathcal L_\ell(X)+Q^{(\ell)}=0$ as in Remark~\ref{rem:zero-terminal-input}.

            \item If $n_\ell>0$ and $m_\ell>0$, solve the stochastic Lyapunov
      equation \eqref{eq:initial-feedback-lyapunov} where $\Phi=\Phi^{(\ell)}$ is the stabilizing terminal feedback obtained from the deflation recursion and apply the
      Newton--Lyapunov iteration \eqref{eq:newton-lyapunov}
      until convergence to $X_-^{(\ell)}$.
        \end{itemize}

\item Factor $X_-^{(\ell)}=Z_\ell Z_\ell^\top$ and apply the factor recursion
      \eqref{eq:factor-reconstruction} for $k=\ell-1,\ldots,0$. Then
      $X_-^{(0)}=Z_0Z_0^\top$ and
      \[
      H_{\rm a}(x)=\tfrac12\|Z_0^\top x\|^2.
      \]
\end{enumerate}

\begin{theorem}
Assume the hypotheses of Proposition~\ref{prop:available-storage}, including
the observability hypothesis ensuring $X_-\succ0$. Let $\ell$ denote the terminal deflation level and assume
Assumption~\ref{ass:regular-feasibility} if $n_\ell>0$ and $m_\ell>0$.
Then the deflation reaches level $\ell$ after finitely many steps,
the terminal computation yields the minimal terminal KYP solution,
and the factor recursion reconstructs the available-storage
matrix of the original system.
\end{theorem}

\begin{proof}
Finite termination and exact transport of the minimal solution follow from
Theorem~\ref{thm:deflation-step} and \eqref{eq:order-equivalence}. Stability is
preserved by Theorem~\ref{thm:stabilizer-transport}.
For $n_\ell>0$ and $m_\ell>0$, the terminal Riccati step is covered by
Theorem~\ref{thm:riccati-storage-identification}.
The remaining terminal cases are covered by
Remarks~\ref{rem:zero-terminal-input} and~\ref{rem:zero-terminal-state}.
In all three cases, the factor recursion
\eqref{eq:factor-reconstruction} is equivalent to
\eqref{eq:minimal-reconstruction}.
\end{proof}

\section{Low-rank numerical realization}
\label{sec:numerics}

For large sparse systems we represent the Newton iterates and the matrices
created by deflation by thin symmetric factorizations
\[
X\approx LD_XL^\top,
\]
with thin $L$ and small symmetric $D_X$; after compression $D_X$ may be taken
diagonal with entries $\pm1$. This signed representation is needed because the
initial and intermediate Newton matrices need not be positive semidefinite.
Only where positivity is known do we use $ZZ^\top$ factors. Generalized
Lyapunov equations and their low-rank treatment are discussed in
\cite{BennerBreiten2013,BennerDamm2011,BreitenRingh2019,ShankSimonciniSzyld2016};
for indefinite right-hand sides see also \cite{LangMenaSaak2015}.

Following the projection principle for deterministic large-scale KYP
deflation \cite{PoloniReis2012}, we exploit that the loop transformation
modifies the state operators only through low-rank terms:
$\widehat A=A+BF$ and $\widehat A_j=A_j+B_jF$.

\subsection{Factorized deflation}

At deflation level $k$, let
\[
H_k=(E^{(k)})^\top J^{(k)}\succeq0.
\]
If
\[
H_k=U_k
\begin{bmatrix}\Lambda_k&0\\0&0\end{bmatrix}U_k^\top,
\qquad \Lambda_k\succ0,
\]
and $U_{k,1}$ contains the eigenvectors belonging to $\Lambda_k$, then
\[
X_{\rm p}^{(k)}=Z_{{\rm p},k}Z_{{\rm p},k}^\top,
\qquad
Z_{{\rm p},k}=J^{(k)}U_{k,1}\Lambda_k^{-1/2}.
\]
Thus $Z_{{\rm p},k}$ has at most $\rank E^{(k)}$ columns.

The remaining reduced coefficients can also be generated from thin factors.
Suppressing the deflation index, assume
\[
Q=C_QD_QC_Q^\top,
\qquad
X_{\rm p}=Z_{\rm p}Z_{\rm p}^\top.
\]
Then
\[
\Theta=C_\Theta D_\Theta C_\Theta^\top,
\]
where
\begin{align*}
C_\Theta=&
\begin{bmatrix}
C_Q&A^\top Z_{\rm p}&Z_{\rm p}&
A_1^\top Z_{\rm p}&\cdots&A_q^\top Z_{\rm p}
\end{bmatrix},\\
D_\Theta=&
\operatorname{diag}\!\left(
D_Q,
\begin{bmatrix}0&I\\ I&0\end{bmatrix},
I,\ldots,I
\right),\\
\Xi=&
Z_{\rm p}(Z_{\rm p}^\top BU_1)
+\sum_{j=1}^q(A_j^\top Z_{\rm p})(Z_{\rm p}^\top B_{j1})+S_1.
\end{align*}
Hence the data in \eqref{eq:reduced-data} require only sparse
matrix--thin-matrix products and small dense algebra.

The stabilizer transport admits an ambient-space representation in terms of the
projector onto the reduced state space.
With the orthonormal basis $\bar V$ of $\im E$ introduced above, set
\[
\Pi=I-\bar V\bar V^\top,
\qquad M_P=\bar V^\top P\bar V,
\]
and let $\Phi_1=U_1^\top\Phi$. The lifted reduced feedback and certificate are
\[
\widetilde\Phi^+=
\begin{bmatrix}
-M_P^{-1}\bar V^\top P\Pi\\[1mm]
\Phi_1(I-\bar V M_P^{-1}\bar V^\top P)\Pi
\end{bmatrix},\qquad
\widetilde P^+=
\Pi(P-P\bar V M_P^{-1}\bar V^\top P)\Pi.
\]
These equal $\Phi^+V^\top$ and $VP^+V^\top$, respectively. Hence the transport step requires only $P\bar V$ and $\bar V^\top P\bar V$ in addition to
thin subspace factors. In the benchmark below the certificate is explicit,
$P=I$, so these quantities are available directly.

\subsection{Sparse realization of the extracted Riccati equation}
\label{subsec:sparse-riccati}

With $\bar V$ as above, set
\[
\Pi=I-\bar V\bar V^\top=VV^\top.
\]
Represent a reduced matrix $Z$ in the ambient state space by
\[
\widetilde X=VZV^\top,
\qquad \widetilde X=\Pi\widetilde X\Pi,
\]
and use the ambient-space representatives
\begin{align*}
\widetilde A=&\Pi A\Pi,
&\widetilde A_j=&\Pi A_j\Pi,&\widetilde Q=&\Pi\Theta\Pi,\\
\widetilde B=&
\begin{bmatrix}\Pi A\bar V&\Pi BU_1\end{bmatrix},&
\widetilde B_j=&
\begin{bmatrix}\Pi A_j\bar V&\Pi B_{j1}\end{bmatrix},
&\widetilde S=&
\begin{bmatrix}\Pi\Theta\bar V&\Pi\Xi\end{bmatrix},
&\widetilde R=&R^+.
\end{align*}
They satisfy, for example,
$\widetilde A=VA^+V^\top$, $\widetilde B=VB^+$, and
$\widetilde Q=VQ^+V^\top$, so this is exactly the coordinate-reduced problem
represented on $\im\Pi$.

At the terminal level, $\widetilde R\succ0$. Define
\[
\widetilde{\mathcal G}(\widetilde X)
=
\widetilde X\widetilde B
+\sum_{j=1}^q\widetilde A_j^\top\widetilde X\widetilde B_j
+\widetilde S,
\qquad\widetilde R_{\widetilde X}
=
\widetilde R-
\sum_{j=1}^q\widetilde B_j^\top\widetilde X\widetilde B_j.
\]
On $\widetilde R_{\widetilde X}\succ0$, the extracted terminal Riccati equation
is
\[
0={}
\widetilde A^\top\widetilde X+\widetilde X\widetilde A
+\sum_{j=1}^q\widetilde A_j^\top\widetilde X\widetilde A_j
+\widetilde Q+\widetilde{\mathcal G}(\widetilde X)
\widetilde R_{\widetilde X}^{-1}
\widetilde{\mathcal G}(\widetilde X)^\top,
\qquad
\widetilde X=\Pi\widetilde X\Pi.
\]
The projection-based formulation extends recursively to successive
deflation steps. The sparse implementation used in the numerical
experiments is restricted to a single deflation step for which
$\im E$ is a coordinate subspace.

For a lifted feedback $\Psi=\Psi\Pi$, set
\[
\widetilde Q_\Psi
=
\widetilde Q+\widetilde S\Psi+\Psi^\top\widetilde S^\top
-\Psi^\top\widetilde R\Psi.
\]
The transported terminal feedback determines the initial stochastic Lyapunov equation.
For a Newton iterate $\widetilde X_i$, define
\[
\widetilde\Gamma_i=
\widetilde R_{\widetilde X_i}^{-1}
\widetilde{\mathcal G}(\widetilde X_i)^\top,
\qquad
\widetilde A_{{\rm cl},i}=\widetilde A+\widetilde B\widetilde\Gamma_i,
\qquad\widetilde N_{j,i}=\widetilde A_j+\widetilde B_j\widetilde\Gamma_i.
\]
Then Newton's equation is
\begin{equation}
\label{eq:projected-newton-lyap}
\widetilde A_{{\rm cl},i}^\top\widetilde X_{i+1}
+\widetilde X_{i+1}\widetilde A_{{\rm cl},i}
+\sum_{j=1}^q\widetilde N_{j,i}^\top
\widetilde X_{i+1}\widetilde N_{j,i}
+\widetilde Q_{\widetilde\Gamma_i}=0,
\qquad
\widetilde X_{i+1}=\Pi\widetilde X_{i+1}\Pi.
\end{equation}
Thus the Riccati iteration reduces to generalized Lyapunov solves on the
terminal regular subspace.

We store
\[
\widetilde X_i\approx L_iD_i^XL_i^\top,
\qquad \Pi L_i=L_i,
\]
and form all products through $L_i$. If
$\widetilde Q=C_QD_QC_Q^\top$, then
\[
\widetilde Q_{\widetilde\Gamma_i}=C_iD_i^QC_i^\top,
\qquad
C_i=\begin{bmatrix}C_Q&\widetilde S&\widetilde\Gamma_i^\top\end{bmatrix},
\qquad D_i^Q=
\begin{bmatrix}
D_Q&0&0\\
0&0&I\\
0&I&-\widetilde R
\end{bmatrix}.
\]
The width of this factor is therefore controlled by the rank of
$\widetilde Q$ and the terminal input dimension.

For the experiments, \eqref{eq:projected-newton-lyap} is solved by the
fixed-point iteration
\[
\widetilde A_{{\rm cl},i}^\top Y_{\nu+1}+Y_{\nu+1}\widetilde A_{{\rm cl},i}
+\widetilde Q_{\widetilde\Gamma_i}
+\sum_{j=1}^q\widetilde N_{j,i}^\top Y_\nu\widetilde N_{j,i}=0,
\qquad Y_{\nu+1}=\Pi Y_{\nu+1}\Pi,\quad Y_0=0.
\]
If $Y_\nu=L_\nu D_\nu L_\nu^\top$, each noise term has the thin factorization
\[
\widetilde N_{j,i}^\top Y_\nu\widetilde N_{j,i}
=(\widetilde N_{j,i}^\top L_\nu)D_\nu
(\widetilde N_{j,i}^\top L_\nu)^\top.
\]
Every inner step is therefore a standard Lyapunov equation with a thin,
possibly indefinite right-hand side, which we solve by signed low-rank ADI.
The splitting is the projected-space analogue of the
stationary iteration in
\cite[Section~2, Eq.~(2.3) and Algorithm~2.1]{ShankSimonciniSzyld2016}.
For the convergence argument, work on the reduced
state space of dimension $n-r$. Let
\[
A_{{\rm cl},i}^+
=V^\top\widetilde A_{{\rm cl},i}V,
\qquad
N_{j,i}^+
=V^\top\widetilde N_{j,i}V,
\]
and define, for $Y\in\mathbb S^{n-r}$,
\[
\mathcal L_0(Y)
=(A_{{\rm cl},i}^+)^\top Y+YA_{{\rm cl},i}^+,\qquad
\mathcal N(Y)
=\sum_{j=1}^q(N_{j,i}^+)^\top YN_{j,i}^+.
\]
Exponential mean-square stability of the closed-loop system
implies, by \cite[Theorem~3.6.1]{Damm2004}, that
$\mathcal L_0$ is invertible and that
\[
\rho(\mathcal T)<1,
\qquad
\mathcal T=-\mathcal L_0^{-1}\mathcal N.
\]
Consequently, the fixed-point iteration converges for exact
inner Lyapunov solves, independently of the definiteness
of the right-hand side. In the low-rank realization,
the inner equations are solved approximately, with convergence
monitored through the corresponding residuals.

In direct ADI methods for generalized stochastic Lyapunov equations,
multiplicative-noise operators affect convergence, and the classical
rational minimax shift problem does not extend directly. Both
Wachspress shifts and bilinear $\mathcal H_2$ interpolation points
have been investigated in this setting \cite{BennerBreiten2013}.
Our fixed-point splitting instead applies ADI only to standard
Lyapunov equations. Shift selection for standard low-rank ADI and
Newton--ADI methods, including adaptive and residual-minimizing
strategies, is discussed in
\cite{BennerKuerschnerSaak2014,BennerMenaSaak2008,Kuerschner2019}.
Here we use a fixed, heuristically chosen set of real negative shifts
throughout the Newton and fixed-point iterations, without explicitly
including the multiplicative-noise operators in shift selection.
This avoids repeated shift generation and permits reuse of sparse
shifted factorizations. The Newton-dependent feedback corrections
are handled by the Sherman--Morrison--Woodbury formula below.

Projected shifted solves can be performed with the sparse full-order operator.
For $b\in\im\Pi$, the system
\[
(\Pi A^\top\Pi+pI)x=b,
\qquad x\in\im\Pi,
\]
is equivalent to
\[
\begin{bmatrix}
A^\top+pI&\bar V\\
\bar V^\top&0
\end{bmatrix}
\begin{bmatrix}x\\\lambda\end{bmatrix}
=
\begin{bmatrix}b\\0\end{bmatrix}.
\]
The large block retains the sparsity of $A$, while $\bar V$ is only a thin
augmentation. The closed-loop shift is a low-rank modification,
\[
\widetilde A_{{\rm cl},i}^\top+pI
=
\widetilde A^\top+pI+\widetilde\Gamma_i^\top\widetilde B^\top.
\]
Let $p<0$ be a real ADI shift such that
$-p\notin\sigma(A^+)$.
Then the projected open-loop shift
$(\widetilde A^\top+pI)|_{\im\Pi}$ is invertible;
denote its inverse by $\mathcal S_p$.
Since the closed-loop drift on $\im\Pi$ is Hurwitz,
the corresponding closed-loop shift is also invertible.
Applying the standard Sherman--Morrison--Woodbury identity
to these operators on $\im\Pi$ gives
\[
\begin{aligned}
\Bigl(\left.(\widetilde A_{{\rm cl},i}^\top+pI)\right|_{\im\Pi}\Bigr)^{-1}
={}&\mathcal S_p
-\mathcal S_p\widetilde\Gamma_i^\top
\bigl(I+\widetilde B^\top\mathcal S_p\widetilde\Gamma_i^\top\bigr)^{-1}
\widetilde B^\top\mathcal S_p .
\end{aligned}
\]
The additional dense inverse has only the terminal input dimension.

\subsection{Compression, residuals, and reconstruction}\label{subsec:compression-residuals}

A symmetric factor $C_{\rm f}D_{\rm f}C_{\rm f}^\top$ is compressed by a thin
QR factorization $C_{\rm f}=Q_{\rm f}R_{\rm f}$ followed by an eigendecomposition of
the small matrix $R_{\rm f}D_{\rm f}R_{\rm f}^\top$. For residual evaluation, each Lyapunov, Riccati or KYP
residual is assembled as a signed factorization
$C_{\rm res}D_{\rm res}C_{\rm res}^\top$,
without forming the dense residual matrix.
The same QR reduction yields its Frobenius norm and
nonzero eigenvalues from the small symmetric matrix
$R_{\rm res}D_{\rm res}R_{\rm res}^\top$.
If the residual factor is rank deficient in the ambient
space, zero must also be included when determining
an extreme eigenvalue.

After convergence, positivity permits
\[
X_-^{(\ell)}\approx Z_\ell Z_\ell^\top.
\]
Backward reconstruction is performed factorwise using
\eqref{eq:factor-reconstruction}. In the projected realization, the second block is kept already lifted to the
current ambient space, so no dense $V^{(k)}$ is formed; the concatenated factor
is recompressed after each step.

For the chain benchmark below,
\[
E=[e_1\ 0],
\qquad \bar V=e_1,
\qquad \Pi=\operatorname{diag}(0,1,\ldots,1),
\]
so projection deletes the first state coordinate. The terminal matrices
$A^+$ and $A_1^+$ in the chain benchmark remain banded and are stored sparsely; the ADI
solves use reusable sparse factorizations of $(A^+)^\top+pI$.

Under observability the exact available-storage matrix is positive definite
and therefore generally full rank. Here ``low-rank'' refers to the compressed
numerical factors used throughout the computation; their effectiveness depends
on numerical rather than algebraic rank. Relevant decay results for generalized
Lyapunov equations are discussed in \cite{BennerBreiten2013}.

\section{Numerical experiments}
\label{sec:experiments}

We first verify the low-rank realization against a dense reference computation
and then study its behavior under increasing state dimension. The example is
chosen so that passivity, singularity of the feedthrough block, one nontrivial
deflation step, and terminal regular feasibility can all be checked
analytically.

\subsection{A sparse stochastic passive chain}

Let $e_i$ denote the $i$th Euclidean unit vector, consider $n\ge 3$, and define the skew-symmetric
nearest-neighbor matrix
\[
\mathsf J_n=
\sum_{i=1}^{n-1}c_i
\bigl(e_ie_{i+1}^\top-e_{i+1}e_i^\top\bigr),
\qquad c_i>0.
\]
We use
\[
c_i=1+0.2\sin i,
\qquad
h=\left\lfloor\frac{n+1}{2}\right\rfloor,
\]
and set
\[
A_1=\sigma e_he_h^\top,
\qquad
A=\mathsf J_n-\delta I-\frac12A_1^\top A_1.
\]
The two input channels are
\[
B=\begin{bmatrix}e_n&e_1\end{bmatrix},
\qquad
B_1=\begin{bmatrix}\gamma e_n&0\end{bmatrix},
\qquad
C=B^\top,
\qquad
D=\frac12
\begin{bmatrix}
\rho&0\\0&0
\end{bmatrix}.
\]
All experiments use
\[
\delta=0.2,
\qquad
\sigma=0.15,
\qquad
\gamma=0.35,
\qquad
\rho=1,
\qquad
\kappa=0.5.
\]
Since $2\le h\le n-1$, one has $A_1^\top B_1=0$, and therefore $X=I$ is an
explicit KYP solution:
\[
A^\top+A+A_1^\top A_1=-2\delta I,
\qquad
B-C^\top+A_1^\top B_1=0,
\]
while
\[
D+D^\top-B_1^\top B_1
=
\begin{bmatrix}
\rho-\gamma^2&0\\0&0
\end{bmatrix}\succeq0.
\]
The second input direction is thus genuinely singular. Moreover, $(C,A)$ is
observable. Indeed, the second output component already suffices: the first new
nonzero entry of $e_1^\top A^k$ occurs in position $k+1$ and has coefficient
$c_1\cdots c_k\ne0$.

We apply the loop transformation with $K=\kappa I_2$ and put
\[
\eta=\frac{1}{1+\kappa\rho/2}.
\]
The singular transformed input produces the affine constraint $Xe_1=e_1$ and
hence
\[
X_{\rm p}=e_1e_1^\top.
\]
For the loop-transformed zero feedback, $P=I$ is a strict Lyapunov certificate.
Indeed, a direct calculation gives
\[
\widehat A^\top+\widehat A+\widehat A_1^\top\widehat A_1
=
-2\delta I-2\kappa e_1e_1^\top
+\bigl(-2\kappa\eta+\gamma^2\kappa^2\eta^2\bigr)e_ne_n^\top
\prec0.
\]
Since $\bar V=e_1$, the off-diagonal block $P_{12}$ in the transport formula
vanishes. Hence $\Lambda=0$, and the stabilizing terminal feedback produced by
Theorem~\ref{thm:stabilizer-transport} is again the zero feedback. This is the
feedback used for the terminal initialization below.

After this single deflation step, the terminal problem has state dimension
$n-1$ and two input channels. With
$\bar e_1,\bar e_{n-1}\in\R^{n-1}$ denoting the boundary unit vectors,
the reduced input matrices are
\[
B^+=\begin{bmatrix}-c_1\bar e_1&\eta\bar e_{n-1}\end{bmatrix},
\qquad
B_1^+=\begin{bmatrix}0&\gamma\eta\bar e_{n-1}\end{bmatrix}.
\]
Since $V^\top X_{\rm p}=0$ and
$\widehat A_1^\top X_{\rm p}\widehat A_1=0$,
the definition of $\Theta$ yields $Q^+=V^\top\widehat QV$.
Using
\[
\widehat Q
=2\kappa e_1e_1^\top
 +(2\kappa\eta-\kappa^2\rho\eta^2)e_ne_n^\top
\]
and $\eta^{-1}=1+\kappa\rho/2$, we obtain
\[
Q^+=2\kappa\eta^2\bar e_{n-1}\bar e_{n-1}^\top,
\qquad
R^+=
\begin{bmatrix}
2\delta&0\\0&\rho\eta^2
\end{bmatrix}\succ0.
\]
In particular, $Q^+$ has rank one. Since $I_n=X_{\rm p}+VI_{n-1}V^\top$ and
$I_n$ satisfies the original KYP inequality,
the deflation theorem implies that $I_{n-1}$
is feasible for the reduced problem.
It is a regular-feasibility witness because
\[
R^+_{I}=
\begin{bmatrix}
2\delta&0\\
0&(\rho-\gamma^2)\eta^2
\end{bmatrix}\succ0.
\]
Thus all assumptions used in the terminal Riccati theory are verified directly
for the benchmark. The sparse specialization is exactly the coordinate-deletion
case described at the end of Subsection~\ref{subsec:compression-residuals}.

\subsection{Dense verification}

For $n=100$, we first solved the terminal problem by dense arithmetic. The
initial stochastic Lyapunov solve attained a relative generalized-Lyapunov residual of
$7.8\cdot10^{-15}$, and Newton's method required five updates to reach a
terminal Riccati residual of $1.2\cdot10^{-13}$. The corresponding low-rank
implementation produced the same Newton branch with relative Frobenius error
\[
\frac{\|X_{\rm LR}-X_{\rm dense}\|_{\rm F}}
     {\|X_{\rm dense}\|_{\rm F}}
=5.54\cdot10^{-9}
\]
and relative feedback error $3.44\cdot10^{-9}$. After backward reconstruction,
the largest eigenvalue of the original KYP matrix was $5.43\cdot10^{-9}$.

The final dense Newton residuals
\[
6.61\cdot10^{-2},\qquad
6.46\cdot10^{-4},\qquad
7.98\cdot10^{-8}
\]
are consistent with local quadratic convergence until the accuracy
of the inner solves becomes dominant. At the relative compression
threshold $10^{-11}$, the terminal low-rank factor has rank $86$.

\subsection{Scaling and low-rank behavior}

For the scaling study, we used the same parameters and numerical tolerances for
all state dimensions. All computations were performed in double precision. We used twelve cyclic real negative ADI shifts whose absolute
values were logarithmically spaced between $0.2$ and $3$,
with inner ADI tolerance $10^{-9}$. The generalized Lyapunov equations
were terminated at a relative residual of $10^{-8}$, low-rank factors were
compressed with relative threshold $10^{-11}$, and Newton's iteration was
stopped when
\[
\frac{\|\mathcal R(X)\|_{\rm F}}
{\max\{1,\|\mathcal H(X)\|_{\rm F},
\|\mathcal G(X)R_X^{-1}\mathcal G(X)^\top\|_{\rm F}\}}
\le 2\cdot10^{-8}.
\]
Using MATLAB's \texttt{tic}/\texttt{toc}, we timed preprocessing, the
low-rank Riccati--Newton iteration, backward reconstruction, and the
factorized KYP audit, excluding system construction, plotting, and file
output. An unmeasured warm-up was followed by five timed runs at each
dimension; Table~\ref{tab:scaling} reports their median. Its KYP column
gives the largest eigenvalue of the original KYP matrix after factorized
reconstruction, so small positive values quantify the numerical
violation of the KYP inequality.

All timings were obtained in MATLAB R2026a Update 2 on an Intel Core Ultra 7
258V processor with 8 computational threads and 32~GB LPDDR5X-8533 memory.
MATLAB used Intel oneAPI MKL 2025.0.1 for BLAS/LAPACK with AVX2, and
\texttt{maxNumCompThreads} was set to~8. The MATLAB implementation was generated using ChatGPT (OpenAI)
and thoroughly reviewed and validated by the authors.
\begin{table}[ht]
\centering
\caption{Scaling of the low-rank deflation--Riccati iteration. Here $r_0$ is
the compressed rank after the initial stochastic Lyapunov solve and $r_*$ is the
terminal compressed factor rank.}
\label{tab:scaling}

\medskip
\begin{tabular}{r|rrrrrr}
$n$ & $r_0$ & $r_*$ & Newton & Riccati residual & KYP max. eig. & time [s]\\
\hline
100  & 61 & 86 & 5 & $8.10\cdot10^{-9}$ & $5.43\cdot10^{-9}$ & 1.11\\
250  & 44 & 88 & 5 & $7.89\cdot10^{-9}$ & $5.25\cdot10^{-9}$ & 1.66\\
500  & 44 & 88 & 5 & $7.55\cdot10^{-9}$ & $4.83\cdot10^{-9}$ & 2.22\\
1000 & 44 & 87 & 5 & $8.06\cdot10^{-9}$ & $5.32\cdot10^{-9}$ & 3.11\\
2000 & 44 & 87 & 5 & $7.63\cdot10^{-9}$ & $4.93\cdot10^{-9}$ & 4.97\\
5000 & 44 & 88 & 5 & $7.88\cdot10^{-9}$ & $5.22\cdot10^{-9}$ & 14.63
\end{tabular}
\end{table}

From $n=250$ onward, the initial factor rank stays at $44$ and the
terminal rank at $87$ or $88$, despite increasing state dimension.
For $n=5000$, the terminal and reconstructed factors have $88$ and $89$
columns, respectively. All $88$ eigenvalues of the compressed terminal
core are positive, and the effective input matrix remains uniformly
regular, with $\lambda_{\min}(R_X)=0.4$ in every run.

All dimensions require five Newton updates, while the Riccati and KYP
residuals remain nearly constant. Median runtime grows from $1.11$ to
$14.63$ seconds, a factor of about $13.2$ for a fifty-fold increase in
$n$. Although these timings do not establish an asymptotic complexity
law, they demonstrate computations with thousands of states without
dense $n\times n$ matrices: dense square-matrix operations involve only
input- or factor-rank dimensions, while the low-rank factors have $n$ rows.

\section{Conclusions}

We have developed a structure-preserving approach to the computation
of available storage for stochastic passive systems with singular
feedthrough and control-dependent diffusion. Structural deflation
reduces the stochastic KYP inequality without regularization,
preserving its Loewner-minimal solution and transporting stabilizing
feedbacks. Under terminal regular feasibility, the resulting rational
stochastic Riccati equation can be solved by a monotonically convergent
Newton iteration initialized through a single stochastic Lyapunov
equation. Degenerate terminal cases are handled without a Riccati solve.

The sparse low-rank implementation combines projected operators,
factorized deflation, and generalized stochastic Lyapunov solves.
Numerical experiments on a stochastic chain with up to $n=5000$
states show convergence within five Newton updates and
nearly constant compressed factor ranks. In particular, these
results show that the available-storage matrix can be computed
efficiently in low-rank form even when the exact solution is
positive definite and therefore has full rank.

\section*{Acknowledgements}

Funded by the Deutsche Forschungsgemeinschaft (DFG, German Research
Foundation) -- Project-ID 531152215 -- SFB 1701 ``Port-Hamiltonian Systems''.

\section*{CRediT authorship contribution statement}

Julia Ackermann and Timo Reis: Conceptualization, Methodology,
Formal analysis, Software, Validation, Writing -- original draft,
Writing -- review \& editing.

\section*{Declaration of competing interests}

The authors declare that they have no known competing financial
interests or personal relationships that could have appeared to
influence the work reported in this paper.

\section*{Data and code availability}

The MATLAB implementation and scripts for reproducing the numerical
experiments are provided as supplementary material.
No external datasets were used.

\section*{Declaration of generative AI and AI-assisted technologies}

All theoretical results and proofs were developed by the authors.
ChatGPT (OpenAI) was used for occasional language polishing and
for generating the MATLAB code used in the numerical experiments.
The AI-generated code was thoroughly reviewed and validated by the
authors, who take full responsibility for the content and results
of this work.

\bibliographystyle{plain}
\bibliography{references}

\end{document}